\documentclass[reqno]{amsart}
\usepackage[utf8]{inputenc}
\usepackage[T1]{fontenc}
\usepackage{enumerate,amssymb}
\usepackage{lmodern, microtype}
\usepackage{enumerate, amscd, color}

\usepackage{tikz}
\usetikzlibrary{decorations.pathreplacing}
\usepackage{color}

\definecolor{darkgreen}{rgb}{0,0.6,0}
\definecolor{darkred}{rgb}{0.7,0,0}
\definecolor{darkblue}{rgb}{0,.1,.6}
\definecolor{darkgray}{rgb}{0.3,.3,.3}
 \usepackage[%
   colorlinks, 
  citecolor=darkgreen, linkcolor=darkblue,
   pdfpagelabels=true,
   unicode=true,
   pdfusetitle
 ]{hyperref}

\usepackage[normalem]{ulem}
\renewcommand\sout{\bgroup\markoverwith
{\textcolor{red}{\rule[0.7ex]{3pt}{1.4pt}}}\ULon}
 
\newcommand\II{\mathrm{I}\hskip-.3mm\mathrm{I}}

\newcommand\datver[1]{\def\datverp
{\par\boxed{\boxed{\text{Version: #1; Run: \today}}}}}\datver{0.1}

\newcommand{\define}{\mathrel{\mathrm{:=}}}

\newcommand{\NN}{\mathbb N}

\newcommand{\RR}{\mathbb R}

\newcommand{\ZZ}{\mathbb Z}

\newcommand\pa{\partial}

\newcommand{\loc}{\operatorname{loc}}

\newcommand{\vol}{\operatorname{vol}}
\newcommand{\dvol}{\operatorname{dvol}}

\newcommand{\maK}{\mathcal K}

\newcommand{\maV}{\mathcal V}

\newcommand{\rinj}{\mathop{r_{\mathrm{inj}}}}

\newtheorem{theorem}{Theorem}

\newtheorem{lemma}[theorem]{Lemma}

\theoremstyle{definition}
\newtheorem{definition}[theorem]{Definition}
\theoremstyle{remark}
\newtheorem{remark}[theorem]{Remark}
\newtheorem{example}[theorem]{Example}

\author[B. Ammann]{Bernd Ammann} \address{B. Ammann, Fakult\"at f\"ur
  Mathematik, Universit\"at Regensburg, 93040 Regensburg, Germany}
\email{bernd.ammann@mathematik.uni-regensburg.de}

\author[N. Gro{\ss}e]{Nadine Gro{\ss}e} \address{N. Gro{\ss}e,
  Mathematisches Institut, Universit\"at Freiburg, 79104 Freiburg,
  Germany} \email{nadine.grosse@math.uni-freiburg.de}

\author[V. Nistor]{Victor Nistor}\address{V. Nistor, Universit\'{e} de
  Lorraine, CNRS, IECL, F-57000 Metz, France
and Inst. Math. Romanian Acad.  PO BOX 1-764, 014700 Bucharest Romania}
\email{victor.nistor@univ-lorraine.fr}

\thanks{B.A. and N.G. have been partially supported by SPP 2026
  (Geometry at infinity), funded by the DFG.  B.A. has also been
  partially supported by the DFG SFB 1085 (Higher
  Invariants). V.N. has been partially supported by
  ANR-14-CE25-0012-01.}
\begin{document}

\title[Domains with singular points]{Analysis and boundary value problems
  on singular domains: an approach via bounded geometry}

\begin{abstract}
We prove well-posedness and regularity results for elliptic boundary
value problems on certain singular domains that are conformally
equivalent to manifolds with boundary and bounded geometry. Our
assumptions are satisfied by the domains with a smooth set of singular
cuspidal points, and hence our results apply to the class of domains
with isolated oscillating conical singularities.  In particular, our
results generalize the classical $L^2$-well-posedness result of
Kondratiev for the Laplacian on domains with conical points. However,
our domains and coefficients are too general to allow for singular
function expansions of the solutions similar to the ones in
Kondratiev's theory. The proofs are based on conformal changes of
metric, on the differential geometry of manifolds with boundary and
bounded geometry, and on our earlier geometric and analytic results on
such manifolds.
\end{abstract}

\maketitle
\setcounter{page}{1}

\subsection*{Version fran\c{c}aise abr\'eg\'ee}

Nous prouvons des résultats de solvabilité et régularité pour des
systèmes satisfaisants la condition de Legendre forte avec conditions
au bord mixtes de type Dirichlet-Neumann sur certains domaines
singuliers. Notre classe de domaines singuliers contient la classe des
domaines avec des singularités coniques isolées. Nos résultats
généralisent ainsi le théorème d'isomorphisme de Kondratiev
\cite{Kondratiev67}.
Dans la suite, $M$ sera une variété riemannienne lisse à bord de
dimension~$m$ et $ E \to M $ sera un fibré vectoriel hermitien
\'equip\'e d'une connexion. Pour nos résultats, nous allons aussi
supposer que $M$ ait une géométrie bornée. Soit $ a $ une forme
sesquilinéaire lisse sur $T^* M \otimes E $ et $P_a \colon H^1(M; E)
\to H^1(M; E)^*$ défini par la formule $\langle P_au, v \rangle
\define \int_{M} a (\nabla u, \nabla v) d \vol_g $, pour $ u, v \in
H^1 (M; E) $.
Nos espaces de fonctions seront les espaces de Sobolev pondérés de
type Kondratiev, voir Équation \eqref {eq.def.2weights_K}. Nous
supposons donée une partition $\pa M = \pa_0 M \sqcup \pa_1 M $ du
bord en deux sous-ensembles disjoints et ouverts, ainsi que des
conditions au bord $B_j$ d'ordre $j$, $B_j$ sur $\pa_j M$.
Nos résultats sont alors:
\begin{enumerate}[(i)]
\item $ P $ satisfait la régularité dans les espaces pondérés $ f \maK
  _ {(\rho)} ^ {\ell, 2} $ de l'équation ~ \eqref {eq.def.2weights_K}
  si, et seulement si les conditions au bord $B = (B_0, B_1) $
  satisfont la condition de régularité de Shapiro-Lopatinski
  uniforme. Ces conditions sont satisfaites pour les opérateurs
  satisfaisant les conditions de Legendre fortes avec des conditions
  au bord mixtes (Dirichlet/Neumann). On obtient en particulier des
  résultats de régularité pour l'opérateur de Laplace avec conditions
  au bord mixtes, Théorème~\ref{thm.reg}.

\item Si, en plus des conditions de (i), $P$ satisfait une inégalité
  de Hardy-Poincaré, alors le problème au bord associé à $P$ est
  également bien posé.
\end{enumerate}
En principe, la classe des domaines à laquelle nos résultats
s'appliquent est assez large, mais pour des raisons d'espace et afin
de réduire au minimum les détails techniques, nous considérons dans
cette note principalement les exemples de domaines cuspidaux et
wedge. L'ensemble des points singuliers de ces domaines est une
sous-vari\'et\'e lisse compacte.

\subsection*{Introduction} 
We prove well-posedness and regularity results for systems of partial
differential equations satisfying the strong Legendre condition with
mixed Dirichlet-Neumann boundary conditions on certain singular
domains. Our class of singular domains includes the class of domains
with isolated conical singularities and thus they generalize the
classical well-posedness result of Kondratiev
\cite{Kondratiev67}. Unlike Kondratiev's theory, singular functions
expansions are not possible in our setting.

Let us briefly state our main result. Here are first, our
assumptions. Throughout this paper, $(M, g)$ will be a smooth,
$m$-dimensional Riemannian manifold with boundary and $E \to M$ will
be a hermitian vector bundle with connection $\nabla$ such that its
curvature $R^E$ and all its covariant derivatives $\nabla^j R^E$, $j
\ge 1$, are bounded. For our results, we shall also assume that $M$
has bounded geometry (a concept recalled below, see, however
\cite{AGN3, AGN1, BrezisSobolev, GN17} for the concepts not recalled
in this paper). Let $a$ be a smooth, sesquilinear form on $T^*M
\otimes E$ and $P_a \colon H^1(M; E) \to H^1(M; E)^*$ be defined by
the formula $\langle P_au, v \rangle \define \int_{M} a(\nabla u,
\nabla v) d\vol_g$, for $u,v \in H^1(M; E)$. Recall that if $a$ is
uniformly positive definite, then $a$ is said to satisfy the
\emph{strong Legendre condition}. If $P = P_a + Q$, where $Q$ is of
order $\le 1$, we shall say that $P$ satisfies the strong Legendre
condition if, and only if, $a$ does. This implies that $P$ is strongly
elliptic. For scalar operators, the condition that $P$ satisfies the
strong Legendre condition is actually equivalent to $P$ being
uniformly strongly elliptic. A smooth function $f \colon M \to (0,
\infty)$ will be called an \emph{admissible weight} if $f^{-1}df$ has
bounded covariant derivatives of all orders. Let $f$ and $\rho$ be
admissible weights on $M$. If $g_0 \define \rho^2 g$ and $\nabla_0$ is
the Levi-Civita connection associated to $g_0$, we can describe our
function spaces as the following Kondratiev-type weighted Sobolev
spaces
\begin{equation}\label{eq.def.2weights_K}
  f \maK_{(\rho)}^{\ell,p}(M,g_0; E) \define \bigl\{ \psi \ |\ \rho^{j}
  \nabla_0^j ( f^{-1} \psi ) \in L^p(M, g_0 ;T^*M^{\otimes j} \otimes
  E), \ (\forall)\, j \leq \ell\, \bigr\}.
\end{equation}
In our applications and in some of our results, the weight $\rho$ is
bounded. For simplicity, we shall assume this throughout the paper.
We also assume that we have a partition $\pa M = \pa_0 M \sqcup \pa_1
M$ of the boundary in two disjoint, open subsets. We assume that we
are given boundary conditions $B_j$ of order $j$, $B_j$ on $\pa_j M$,
satisfying the boundedness and smoothness conditions stated before
Theorem~\ref{thm.reg}. Our results are then as follows (for $(M, g)$
with bounded geometry):
\begin{enumerate}[(i)]
\item $P$ satisfies regularity in the weighted spaces $f
  \maK_{(\rho)}^{\ell,2}$ of Equation~\eqref{eq.def.2weights_K} if,
  and only if $B = (B_0, B_1)$ satisfies the uniform
  Shapiro-Lopatinski regularity conditions. These conditions are
  satisfied for operators satisfying the strong Legendre conditions
  with mixed (Dirichlet/Neumann) boundary conditions. In particular,
  the Laplace operator satisfies regularity for mixed boundary
  conditions, Theorem~\ref{thm.reg}.
 
\item If, in addition to the conditions of (i), $P$ satisfies a
  Hardy-Poincar\'e inequality, then $P$ also satisfies a
  well-posedness result. We provide several examples of how to obtain
  the Hardy-Poincar\'e inequality.
\end{enumerate}

In principle, the class of domains to which our results apply is
pretty large, but for reasons of space and in order to keep the
technicalities to a minimum, we mostly consider the examples of
canonical cuspidal and wedge domains introduced by H. Amann
\cite{AmannFunctSp}, whose definition is recalled below. The set of
singular points of such a domain is smooth and compact (without
corners).  It is even a finite set for cuspidal domains. Some very
general and nice results were obtained in \cite{DaugeCusp, KMRpoint}
for certain domains with isolated point cusp singularities. Our
methods are quite different, relying heavily on differential geometry,
and thus allowing us to treat a large class of domains. Moreover, our
coefficients are less regular than the ones in the references, but we
lose the Fredholm properties and the singular function expansions
obtained in \cite{DaugeCusp, KMRpoint}. Algebras of pseudodifferential
operators on manifolds with cuspidal points were considered in
\cite{SSS98}. The index problem on such manifolds was considered in
\cite{LeschHorns}. We thank Herbert Amann for useful comments.

\subsection*{Manifolds with boundary and bounded geometry}
In this paper, $(M,g)$ will always be a smooth, $m$-dimensional
Riemannian manifold with boundary and $E \to M$ will be an vector
bundle with metric and metric preserving connection. A smooth function
$f \colon M \to (0, \infty)$ will be called {\em a $g$-admissible
  weight} if $f^{-1}df$ has bounded covariant derivatives of all
orders. We shall say that $E$ has \emph{totally bounded curvature} if
its curvature and all of its covariant derivatives are bounded. We
endow $TM$ with the Levi-Civita connection $\nabla$ associated to
$g$. Recall that $M$ is said to have \emph{bounded geometry} if its
injectivity radius $\rinj(M)>0$ is positive and if $TM$ has
\emph{totally bounded curvature}. We assume from now on that $E$ is
complex and it has totally bounded curvature.

Let us consider a codimension one submanifold $H \subset M$
(i.e. hypersurface). Assume that $H$ carries a globally defined unit
normal vector field $\nu$ and let $\exp^\perp(x,t) \define
\exp^M_x(t\nu_x)$ be the exponential in the direction of the chosen
unit normal vector. By $\II^H$ we denote the \emph{second fundamental
  form} of $H$. The following two definitions are from \cite{AGN1}.

\begin{definition}\label{hyp_bdd_geo} 
Let $(\widehat M, \widehat g)$ be a Riemannian manifold with bounded
geometry. We say that $H \subset \widehat M$ is a \emph{bounded
  geometry hypersurface} in $M$ if it is a closed subset of $M$, if
$\|\nabla^k \II^H \|_{L^\infty} < \infty$ for all $k\geq 0$, and if
here is $r_\pa > 0$ such that $\exp^\perp\colon H\times (-r_\pa, r_\pa
)\to \widehat M$ is a diffeomorphism onto its image.
\end{definition}

\begin{definition}\label{def_bdd_geo}  
A Riemannian manifold~$(M,g)$ with boundary is said to have
\emph{bounded geometry} if there is a Riemannian manifold $\widehat M$
with bounded geometry containing~$M$ as an open subset such that
$\partial M$ is a bounded geometry hypersurface in $\widehat M$.
\end{definition}

\begin{remark}\label{rem.equiv2}
In \cite{AmannFunctSp}, Amann has introduced the class of ``singular
manifolds.'' A \emph{singular manifold} $(M, g_0, \rho)$ is a
Riemannian manifold with boundary $(M, g_0)$ together with a
\emph{singularity function} $\rho$ satisfying suitable properties. In
particular, $(M, \rho^{-2}g_0)$ is assumed to be a manifold with
boundary and bounded geometry. Conversely, if $(M, g)$ has bounded
geometry and $\rho$ is a $g$-admissible weight, then $(M, g_0: =
\rho^2 g)$ is a singular manifold with singularity function $\rho$. In
the boundaryless case, this was first noticed in \cite{Disconzi} (see
also \cite{AmannSingMan}.)  The singularity function $\rho$ is seen to
be a $g$-admissible weight. For manifolds with boundary, this result
follows from \cite{AGN1} or \cite{nadine}. The results of
\cite{AmannFunctSp} apply therefore to the setting of manifolds with
boundary and bounded geometry endowed with an admissible weight. A
triple $(M, g, \rho)$ consisting of a manifold with boundary and
bounded geometry and a bounded $g$-admissible weight $\rho$ will be
called an {\em Amann triple}.
\end{remark}

\subsection*{Conformal changes of metric}
If $h_1, h_2 \colon X \to (0, \infty)$, we shall write $h_1 \sim h_2$
if $h_1/h_2$ and $h_2/h_1$ are both bounded. Let $g_0$ be a second
Riemannian metric on $M$, whose Levi-Civita connection is denoted
$\nabla_0$. Let $\rho, f \colon M \to (0, \infty)$ be measurable
functions and $p\in[1,\infty]$. Recall then from Equation
\eqref{eq.def.2weights_K} the definition of the spaces $f
\maK_{(\rho)}^{\ell,p}(M,g_0; E)$, which reduce to the usual Sobolev
spaces if $\rho, f \sim 1$.
More precisely, if $g \define \rho^{-2} g_0$ and if $\rho$ is a
$g$-admissible weight and $f$ is continuous, then the weighted and
classical spaces are related by
\begin{equation}\label{eq.relation}
  f \maK_{(\rho)}^{\ell,p}(M,g_0; E) = f \rho^{-m/p} W^{\ell,p}(M,g;
  E), \quad 1 \le p \le \infty.
\end{equation}
(see \cite{AmannFunctSp, AmarSobolev, sobolev} and Remark
\ref{rem.equiv2}). We drop the superscript $p$ for $p = 2$:
$\maK_{(\rho)}^{\ell}(M,g_0; E) \define \maK_{(\rho)}^{\ell, 2}(M,g_0;
E)$ and so on. We assume from now on that $g = \rho^{-2}g_0$.

\begin{example}\label{ex.r.lambda}
A typical example is when $M \subset \RR^m$ is the closed unit ball,
$g_0$ is the euclidean metric, and $\rho = r^\lambda$, where $r$ is
the distance to the origin. Then $(M, g := \rho^{-2} g_0)$ has bounded
geometry if, and only if, $\lambda \geq 1$. Moreover,
\begin{equation*} f \ := \
  \begin{cases}
    \ e^{- ( \frac{r}{\epsilon})^{-\epsilon} },
    & \mbox{ if } \lambda = 1+ \epsilon > 1 ,\\
    \ \ r = \rho , & \mbox{ if } \lambda = 1.
  \end{cases}
\end{equation*}
is a $g$-admissible weight. This example is adapted to a domain with
conical points (for instance, a polygonal domain) with set of vertices
$\maV$ by taking $f(x) = \rho(x) := \prod_{P \in \maV} |x - P|$ and
$\lambda = 1$. The extra weight $f$ becomes then unnecessary (for
$\lambda = 1$) and the weighted Sobolev spaces $\maK_{(\rho)}^\ell(M)
\define \{ u \ |\ \rho^{|\alpha|} \pa^\alpha u \in L^{2}(M),
\ (\forall)\, |\alpha| \leq \ell\, \}$ are the spaces introduced by
Kondratiev \cite{Kondratiev67}.
\end{example}

  
We shall \emph{assume from now that $(M, g, \rho)$ is an Amann triple
  (see Remark \ref{rem.equiv2}) and that $f \colon M \to (0, \infty)$
  be a second $g$-admissible weight.} In particular, $\rho$ is a
bounded $g$-admissible weight. We have seen in
Equation~\eqref{eq.relation} how the Sobolev spaces change with
respect to conformal changes of metric. Recall that $g_0 = \rho^2 g$.
For differential operators, a simple calculation based on $L^\infty(M;
E \otimes TM^{\otimes p} \otimes TM^{*\otimes q}, g) =
\rho^{p-q}L^\infty(M; E \otimes TM^{\otimes p} \otimes TM^{*\otimes
  q}, g_0)$ and the fact that $\rho$ is bounded gives:

\begin{lemma}\label{lemma.conf.diff}
Let $P$ be an order $k$ differential operator on $M$ and $P_1 \define
\rho^{k} P$. We have that $P$ satisfies the strong Legendre condition
with respect to the metric $g_0$ if, and only if, $P_1$ satisfies the
strong Legendre condition with respect to the metric $g =
\rho^{-2}g_0$. If $P$ has coefficients in $W^{\infty, \infty}(g_0)$,
then $P_1$ has coefficients in $W^{\infty, \infty}(g)$.
\end{lemma}

A similar result is valid for the boundary differential
  operators appearing as boundary conditions.

\subsection*{Regularity and well-posedness}
Let $P$ be a second order differential operator. We assume from now on
that we are given a partition of the boundary $\pa M = \pa_0 M \sqcup
\pa_1 M$ into two \emph{disjoint, open subsets,} as in \cite{AGN1},
and order $i$ differential boundary conditions $B_i$ on $\pa_i M$.
See, for example, \cite{BrezisSobolev, LionsMagenes1} for general
results on boundary value problems on smooth domains, \cite{DaugeBook,
  KMRpoint, NP, NazarovPopoff} for the case of non-smooth domains, and
\cite{AGN3, GN17} for more general boundary conditions involving
projections. We assume that $\rho^2 P$, $\rho B_1$, and $B_0$
\emph{have coefficients in} $W^{\infty, \infty}(g)$. The typical
assumption is that $P$, $B_1$, and $B_0$ have coefficients in
$W^{\infty, \infty}(g_0)$, which means that they ``stabilize'' towards
the singular points, as in \cite{KMRpoint}, and this is a necessary
condition for the existence of singular function expansions. In view
of Lemma \ref{lemma.conf.diff}, our assumptions are thus weaker, but
singular functions expansions are no longer available in general in
our setting. Our more general setting may be needed in applications to
non-linear PDEs and uncertainty quantification. Also, recall from
\cite{GN17} the uniform Shapiro-Lopatinski regularity conditions and
that they are invariant with respect to conformal changes of
metric. Combining this property with Equation~\eqref{eq.relation} and
with Lemma~\ref{lemma.conf.diff}, we get:

\begin{theorem}\label{thm.reg}
Let $P$ be a $g_0$-uniformly elliptic second order differential
operator acting on sections of $E \to M$ and $B = (B_0, B_1)$ be a
boundary differential operator. We assume that $P$ and $B$ satisfy the
$g_0$-uniform Shapiro-Lopatinski regularity conditions. Then, for any
$\ell \in \NN$, there exists $C > 0$ such that, for all $u \in f
\maK_{(\rho)}^{1}(M,g_0;E)$
\begin{multline*}
  \|u\|_{f \maK_{(\rho)}^{\ell +1}(M,g_0;E)} \leq C \big( \,
  \|Pu\|_{f\rho^{-2} \maK_{(\rho)}^{\ell -1}(M,g_0;E)} + \|u\|_{f
    \maK_{(\rho)}^{1}(M,g_0;E)}\\
  + \|B_0 u\|_{f \rho^{-\frac{1}{2}}\maK_{(\rho)}^{\ell +
      \frac{1}{2}}(\partial_0 M,g_0;E)}
  + \|B_1 u\|_{f \rho^{-\frac{3}{2}}\maK_{(\rho)}^{\ell -
      \frac{1}{2}}(\partial_1 M,g_0;E)} \big).
\end{multline*}
In particular, we can take $P$ to be a uniformly strongly elliptic
scalar operator (such as the Laplacian $P = \Delta_{g_0}$), $B_0u =
u\vert_{\partial_D M}$ (the restriction) and $B_1 u = \pa_\nu^a u$.
\end{theorem}

Let $(M, g_0)$ be a Riemannian manifold with boundary. We now turn to
the well-posedness on $(M, g_0)$. Let $h \colon M \to (0, \infty)$,
and $A \subset \pa M$ be a measurable subset. We shall say that $(M,
A, E, g_0, h)$ satisfies the \emph{Hardy-Poincar\'e inequality} if
there exists a constant $C > 0$ such that, for any $u \in
H^1_{\loc}(M,g_0; E)$, $u = 0$ in $L^2(A)$, we have $\int_{M} |d
u|_{g_0}^2 \dvol_{g_0} \geq C \int_{M} h^{-2} u^2 \dvol_{g_0}.$ The
Hardy-Poincar\'e inequality implies coercive estimates, and hence
well-posedness also for the associated parabolic and hyperbolic
equations, as in \cite{LionsMagenes1}. The Hardy-Poincar\'e inequality
is related to the Poincar\'e inequality, and hence to the concept of
``finite width.''  If $A \subset \pa M$, recall that $(M, A)$ is said
to have \emph{finite width} if $\text{dist}(x, A)$ is uniformly
bounded on $M$ \cite{AGN1} (the distance between two disjoint
connected components of $M$ is $+\infty$). Typically in our results,
the set $A$ will be an open and closed subset of $\pa M$.

\begin{example}\label{ex.r.lambda2}
Again, a typical application is when $g_0$ is the euclidean metric on
$\RR^m$, $r$ is the distance to the origin, and $\lambda > 0$, as in
Example \ref{ex.r.lambda}. However, in this case, we let $\rho =
r^\lambda$ \emph{only} for $r < 1/2$, but set $\rho = r$ for $r > 1$
and $M \subset \RR^m$ is a closed, infinite cone with base a smooth
domain of the unit sphere and with vertex at the origin. Then again,
$(M, g)$ has bounded geometry if, and only if, $\lambda \geq 1$. Also,
$(M, \pa M, g)$ has finite width if, and only if $\lambda =
1$. Finally, $(M, \pa M, g_0, \rho)$ satisfies the Hardy-Poincar\'e
inequality (for $\rho$) if, and only if, $\lambda \leq 1$.
\end{example}

Recall that we have assumed $(M,g)$ to be a Riemannian manifold with
boundary and bounded geometry, $\rho, f \colon M \to (0, \infty)$ to
be $g$-admissible weights, and $g_0 \define \rho^2 g$. We define $P_a$
by $(P_a u, v)_{g_0} \define \int_M a(\nabla u, \nabla v)d\vol_{g_0}$,
with a sesquilinear form $a$ satisfying the \emph{strong Legendre
  condition with respect to $g_0$}. Let $\partial^a_\nu$ be the
conormal derivative associated to $P$, see \cite{GN17}. Combining
Theorem \ref{thm.reg} with the Lax-Milgram Lemma and with the fact
that the Dirichlet and Neumann boundary conditions satify the uniform
Shapiro-Lopatinski regularity conditions \cite{AGN1, GN17}, we obtain:

\begin{theorem}\label{thm.wp}
We assume that $(M, \pa_0 M, E, g_0, \rho)$ satisfies the
Hardy-Poincar\'e inequality. Let $P = P_a$ satisfy the strong Legendre
condition with all $\nabla^j a$ bounded. Then there exists
$\eta_{a,f}>0$ such that, for $|s| < \eta_{a,f}$ and $\ell \geq 1$, we
have an isomorphism
\begin{equation*}
  P_a \colon \rho f^s \maK_{(\rho)}^{\ell + 1}(M, g_0; E) \cap \{ u
  \vert_{\pa_0 M} = 0 , \pa_\nu^a u\vert_{\pa_1 M} = 0 \} \to
  \rho^{-1} f^s \maK_{(\rho)}^{\ell - 1}(M, g_0; E).
\end{equation*}
In particular, we can take $P = \Delta_{g_0}$, the Laplacian
associated to $g_0$. For $\ell = 0$ the result remains true, once one
reformulates it in a variational (i.e. weak) sense.
\end{theorem}

\subsection*{Examples} We include some basic examples.

\subsubsection*{Two dimensional domains}
We consider a (disjoint) partition of the boundary $\pa M = \pa_0 M
\sqcup \pa_1 M$ as above (so $\pa_0 M$ and $\pa_1 M$ are open and
closed).  Recall that $P_a$ is a second order differential operators
on $E \to M$ with coefficients in $W^{\infty, \infty} (g)$ and
satisfying the strong Legendre condition with respect to $g_0$.  For
dimension two domains $M$, the Poincar\'e inequality for $(M, A, g)$
is equivalent to the Poincar\'e inequality for $(M, A, g_0)$ (same
proof as the conformal invariance of the Laplacian in two
dimensions). The Poincar\'e inequality of \cite{AGN1} then gives:

\begin{theorem}\label{thm.wp.2d}
Assume that $(M, \pa_0M, g)$ has finite width and $m \define \dim(M) =
2$. Let $P = P_a$ satisfy the strong Legendre condition with all
$\nabla^j a$ bounded. Then there exists $\eta_{a,f} > 0$ such that,
for $|s| < \eta_{a,f}$ and $\ell \in \ZZ_+$, we have an isomorphism
\begin{equation*}
   P_a \colon \rho f^s \maK_{(\rho)}^{\ell + 1}(M, g_0; E) \cap \{ u
   \vert_{\pa_0 M} = 0 , \pa_\nu^a u\vert_{\pa_1 M} = 0\} \to
   \rho^{-1} f^s \maK_{(\rho)}^{\ell - 1}(M, g_0; E).
\end{equation*}
In particular, we can take $P_a = \Delta_{g_0}$.
\end{theorem}

\subsubsection*{Canonical cuspidal and wedge domains}
We continue with some concrete examples. The simplest examples in
higher dimensions are those of ``model cuspidal and wedge domains.''
We follow the presentation in \cite{AmannFunctSp}. Let $1 < \alpha <
\infty$ and $B \subset \RR^{m-1}$ a compact submanifold, possibly with
boundary, and
\begin{equation}\label{eq.def.cusp}
   K_\alpha^m(B) \define \{ (r, r^\alpha y) \in \RR^m \, \vert
   \ 0 < r < 1, y \in B \},
\end{equation} 
which will be called a \emph{model canonical cusp of order $\alpha$.}
For $\alpha = 1$, we take $B$ a subset of the unit sphere. A domain
with \emph{canonical cuspidal singularities} is a bounded domain
$\Omega \subset \widehat M$ in a Riemannian manifold $(\widehat M,
g_0)$ such that, around each singular point $P$ of the boundary, it is
locally diffeomorphic to $K_{\alpha_P}^m(B_P)$ via a diffeomorphism
defined in a neighborhood of the ambient manifold. Let $\maV$ be the
set of singular points of the boundary, then $\maV$ is finite and we
let $M := \overline{\Omega} \smallsetminus \maV$.  If $\alpha_P = 1$
for all $P \in \maV$, we obtain a domain with \emph{conical points}.
If we replace $K_\alpha^m(B)$ with $K_\alpha^{m-k}(B) \times [0,1]^k$,
$k \ge 0$, we obtain domains with \emph{canonical wedge
  singularities,} in which case, of course, the set $\maV$ of singular
points of $\pa \Omega$ will no longer be finite.

Let us fix $\lambda_P \ge 1$ for each singular point $P \in \maV$. The
weight functions $\rho$ and $f$ are then chosen, around each $P \in
\maV$ as in Example \ref{ex.r.lambda} for $\lambda = \lambda_P$. Let
$g := \rho^{-2} g_0$, as before. If $\lambda_P \geq \alpha_P$ for all
$P$, then $(M, g)$ has bounded geometry (proved in \cite{AmannSingMan}
if $\lambda_P = \alpha_P$ for all $P \in \maV$) and consequently, we
have regularity in the weighted spaces for the mixed Dirichlet-Neumann
problem for operators satisfying the strong Legendre condition. If
$\lambda_P \le \alpha_P$ and $\pa_0 M$ intersects each $V_P$, then
$(M, \pa_0 M, g)$ satisfies the Hardy-Poincar\'e inequality. This
follows from the usual Poincar\'e inequality on each $\{r\} \times
r^{\alpha} B$ by also rescaling in $r$. We work in generalized
spherical coordinates $(r, y) \in (0, \infty) \times S^{m-1}$, so $dx
= r^{m-1} dr dy$. This gives $ r^{-2\alpha} \int_{r^{\alpha} B} |u(r,
y)|^2 dy \le C \int_{r^{\alpha} B} |\nabla_y u(r, y)|^2 dy$ the
Hardy-Poincar\'e inequality $\int_{K_\alpha^m(B)} r^{-2\alpha} |u(r,
y)|^2 dx \le C \int_{K_\alpha^m(B)} |\nabla u(r, y)|^2 dx $ for $u= 0$
on $\pa_1 B$. Let us fix $\lambda_P = \alpha_P$. By considering the
Amann triple $(M, g:= \rho^{-2}g, \rho)$, we obtain that our domain
with canonical wedge singularities satisfies the conclusion
(isomorphism) of Theorem \ref{thm.wp.2d}. For canonical cuspidal
domains and constant coefficient operators, this theorem was first
proved in \cite{KMRpoint}. See also~\cite{NazarovSokolowskiAA09,
  CDN12, DaugeCusp, KamotskiMazya11, MunnierCusp}.

\subsubsection*{Other examples}
Certain simple examples are not ``canonical.'' 

\begin{example}
Let $\Omega = \{(x, y) \in \RR^{2}\ |\ x, y > 0, (y-1)^2 + x^2 >1
\}\cap [-2,2]\times [0,2]$ (see the picture).  The corners
$\{A_i\}_{i=1}^4$ of $\Omega$ are conical points and are treated as
above.
\medskip

\noindent
\begin{minipage}{0.65\textwidth}
Close to $O$ we have two cuspidal open sets ($U$ and its mirror
image), similar to the ones treated in the previous subsection, but
\emph{not canonical}. We have $\lambda = 2$ for these open sets. (We
double $O$, in a certain sense.) Let $r$ be the distance to $O$. Close
to $O$, we then choose $\rho \sim r^{2}$, and $f \sim
e^{-r^{-1}}$. More precisely, $\rho(x) = r^2 \prod_{i=1}^4 |x - A_i|$
and $f = e^{-r^{-1}} \prod_{i=1}^4 |x - A_i|$. Note that we could have
also chosen $\rho \sim x^2$ near $O$.
\end{minipage} \ 
 \hfill \ 
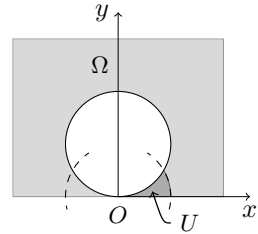
\begin{minipage}{0.3\textwidth}
\centering
\begin{tikzpicture}[scale=0.35]
\draw[fill=gray, opacity=0.3] (-4,-2)--(4,-2) --(4,4) --(-4,4)
--(-4,-2);
\draw (0,3) node[left] {\small $\Omega$};
\draw (0,-2.7) node {\small $O$};
\begin{scope} 
  \clip (0,-2)--(4,-2) --(4,4) --(0,4) --(0 ,-2);
  \draw[fill=gray, opacity=0.5] (0,-2) circle (2cm);
\end{scope}
\draw[line width=0.5, ->] (2,-3) node[right] {\small $U$} .. controls
(1.7,-3) .. (1.3,-1.8);
\draw[fill=white] (0,0) circle (2cm);
\begin{scope} 
  \clip (0,-2) --(4,-3) --(4,4)--(0,-2)--(-4,4)--(-4,-3)--(0,-2);
\draw[dashed] (0,-2) circle (2cm);
\end{scope}
\draw[line width=0.2, ->] (0,-2) --(5,-2) node[below] {$x$};
\draw[line width=0.2, ->] (0,-2) --(0,5) node[left] {$y$};
\end{tikzpicture}
\end{minipage}
\smallskip 

\noindent Let $g \define \rho^{-2}g_E$, where $g_E$ is the standard
(flat) Euclidean metric.  We then have that $(M, g)$, $M \define
\overline{\Omega} \smallsetminus \{O, A_i,i=1,4\}$, is a manifold with
boundary and bounded geometry. Assume that $\pa_0 M$ touches each
singular point (where $O$ is considered as a double point as
above). We can then prove that $(M, \pa_0 M, g_E)$ satisfies the
Hardy-Poincar\'e inequality as in the previous example, and hence
Theorem~\ref{thm.wp.2d} applies. See also \cite{KamotskiMazya11,
  MunnierCusp}.
\end{example}

\begin{example}
Let $f_0, f_1 \colon \RR \to (0, 2\pi)$ satisfy $\|f_i^{(k)}\|_\infty
< \infty$, $k \geq 0$, and $f_1 - f_0 \geq \epsilon > 0$. Let
$\Omega_{f_0, f_1} \define \{ (r \cos \theta, r\sin \theta)\,
\vert\ f_0(\log r) < \theta < f_1(\log r)\}$ and $\Omega$ be a bounded
domain, smooth away from a finite number of points, at which it
coincides, up to a diffeomorphism, with a neighborhood of $0$ in a set
of the form $\Omega_{f_0, f_1}$. Then $\Omega$ is a domain with
\emph{oscillating conical singularities,} similar to the ones studied
by \cite{Rab3}. Theorem~\ref{thm.wp.2d} holds for this domain with
$\rho = f = r = (x^2 + y^2)^{1/2}$. In general, a domain with
\emph{oscillating} conical points \emph{is not} a domain with conical
points.
\end{example}

Finally, in \cite{BMNZ}, it was proved that the assumptions and the
conclusions of Theorem~\ref{thm.wp} are fulfilled by any polyhedral
domain $\Omega \subset \RR^m$ (defined as a suitable stratified space)
with the Euclidean metric $g_0$ for a suitable $g$-admissible weight
$\rho \sim$ the distance to the singular points of the boundary ($g =
\rho^{-2} g_0$).


\def\cprime{$'$}

\end{document}